\documentclass{amsart}

\usepackage{amssymb}
\usepackage{latexsym}
\usepackage{amsmath}
\usepackage{euscript}
\usepackage{graphics}

\newtheorem{theorem}{Theorem}[section]

\newtheorem{lemma}[theorem]{Lemma}
\newtheorem{proposition}[theorem]{Proposition}
\theoremstyle{definition}
\newtheorem{definition}[theorem]{Definition}
\theoremstyle{remark}
\newtheorem{remark}[theorem]{Remark}
\numberwithin{equation}{section}
\renewenvironment {proof} {\begin{trivlist} \item[\hspace{\labelsep}%
\sc Proof.]}{$\Box$ \end{trivlist}}

      \def\dC{{\mathbb C}}
\def\dD{{\mathbb D}}

   \def\dN{{\mathbb N}}   
      \def\dR{{\mathbb R}}
   \def\dT{{\mathbb T}}   
      
   \def\dZ{{\mathbb Z}}

   \def\cH{{\mathcal H}}   
   \def\cK{{\mathcal K}}   
\def\cM{{\mathcal M}}   \def\cN{{\mathcal N}}   
      \def\cR{{\mathcal R}}
\def\cS{{\mathcal S}}   \def\cT{{\mathcal T}}   \def\cU{{\mathcal U}}

\renewcommand {\phi}{\varphi}          

\renewcommand{\Im}{\mathop{\rm Im}\nolimits}        

\def\ran{{\rm ran\,}}                               
\def\ran{{\rm ran\,}}                  

\newcommand {\wt}{\widetilde}          

\def\gh{{\mathfrak h}}

\newcommand{\ptp}{p\times p}
\newcommand{\ptq}{p\times q}
\newcommand{\qtq}{q\times q}

\def\ran{{\rm ran\,}}

\def\ker{{\rm ker\,}}

\def\rank{{\rm rank\,}}

\begin{document}

\begin{title}
{Generalized $\gamma$-generating matrices and Nehari-Takagi problem }
\end{title}
\author{Volodymyr Derkach}
\author{Olena Sukhorukova}

\address{
Kiev, 
Ukraine}
\email{derkach.v@gmail.com}
\email{alena.dn.ua@rambler.ru}

\dedicatory{To the Memory of  Leiba Rodman}

\begin{abstract}
Under certain mild assumption,
we establish a one-to-one correspondence between  solutions of the
Nehari-Takagi problem and solutions of some Takagi-Sarason
interpolation problem.  The resolvent matrix of the Nehari-Takagi
problem is shown to belong to the class of so-called generalized
$\gamma$-generating matrices, which is introduced and studied in the
paper.
\end{abstract}



\maketitle
\section{Introduction}
For a bounded function $f$ defined on $\mathbb{T}=\{z:|z|=1\}$ let us set
\begin{equation}\label{eq:gf_1}
\gamma_k(f)=\frac{1}{2\pi}\int_\dT e^{ik\theta} f(e^{i\theta})d\theta \quad (k=1,2,...).
\end{equation}
The Nehari problem consists of the following: given a sequence of complex numbers $\gamma_k$ $(k\in \mathbb{N})$ find a function $f\in L_\infty(\mathbb{T})$ such that $\|f\|\leq1$ and
\begin{equation}\label{eq:gf_2}
\gamma_k(f)=\gamma_k, \quad (k=1,2,...).
\end{equation}
By Nehari theorem \cite{Neh57} this problem is solvable if and only if the Hankel matrix $\Gamma=(\gamma_{i+j-1})_{i,j=1}^\infty$ determines a bounded operator in $l_2(\dN)$ with $\|\Gamma\|\leq1$.
The problem~\eqref{eq:gf_2} is called indeterminate if it has infinitely many solutions. A criterion for the Nehari problem to be indeterminate and a full description of the set of its solutions was given in  \cite{AAK68}, \cite{AAK71a}.

In \cite{AAK71a} Adamyan, Arov and Krein considered the following indefinite version of the Nehari problem, so called  Nehari-Takagi problem $\mathbf{NTP}_\kappa (\Gamma)$: Given $\kappa\in \mathbb{N}$ and a sequence $\{\gamma_k\}_{k=1}^\infty$ of complex numbers, find a function $f\in L_\infty(\mathbb{T})$, such that $\|f\|_\infty\leq1 $ and
\[
\mbox{rank }(\Gamma(f)-\Gamma)\leq \kappa.
\]
Here $\Gamma (f)$ is the Hankel matrix $\Gamma(f):=(\gamma_{i+j-1}(f))_{i,j=1}^\infty$. As was shown in \cite{AAK71a}, the problem $\mathbf{NTP}_\kappa(\Gamma)$ is solvable if and only if the total multiplicity $\nu_-(I-\Gamma^*\Gamma)$ of the negative spectrum of the operator $I-\Gamma^*\Gamma$ does not exceed $\kappa$. In the case when the operator $I-\Gamma^*\Gamma$  is invertible and $\nu_-(I-\Gamma^*\Gamma)=\kappa$,  the set of solutions of this problem was parameterized by the formula
\begin{equation}\label{eq:factor}
f(\mu)=(a_{11}(\mu)\varepsilon(\mu)+a_{12}(\mu))(a_{21}(\mu)\varepsilon(\mu)+a_{22}(\mu))^{-1},
\end{equation}
where $\mathfrak{A}(\mu)=(a_{ij}(\mu))_{i,j=1}^2$ is the so-called $\gamma$-generating matrix and the parameter $\varepsilon$ ranges over the Schur class of functions bounded and holomorphic on $\mathbb{D}=\{z: \|z\|<1\}$. In \cite{AAK71a}  applications of the Nehari-Takagi problem to various approximation and interpolation problems were presented.
Matrix and operator versions of Nehari problem were considered in \cite{Pa70} and \cite{Ad73}.
In the rational case matrix Nehari and Nehari-Takagi problems were studied in \cite{BGR90}. A complete exposition of these results can be found also in~\cite{Pel03} and~\cite{ArovD08}.

In the present paper we consider the general matrix Nehari-Takagi problem and show that under some assumptions this problem can be reduced to Takagi-Sarason interpolation problem studied earlier in \cite{DD10}. Using the results from~\cite{DD10}, \cite{DD14} we obtain a description of the set of solutions of the matrix Nehari-Takagi problem in the form \eqref{eq:factor}. 

Connections between the class of generalized $\gamma$-generating
matrices and the class of generalized $j$-inner matrix valued
functions (mvf's) introduced in~\cite{DD09} is established. Using
this connection we present another proof of the formula for the
resolvent matrix $\mathfrak{A}(\mu)$ from~\cite{BGR90} in the case
when the Hankel matrix $\Gamma$ corresponds to a rational mvf. All
the results, except the last section, are presented in unified
notations both for the unit circle $\dT$ and the real line $\dR$.

Now we briefly describe the content of the paper.
In Section 2 some preliminary statements concerning generalized Schur mvf's and generalized $j$-inner mvf's are given.
In Section 3 the Takagi-Sarason interpolation problem is studied.
In Section  4 we introduce and study the class 
of generalized $\gamma$-generating matrices and
establish their connection with the class of generalized $j$-inner
mvf's. In Section 5 we consider Nehari-Takagi problem and under
certain mild assumption,  establish a one-to-one correspondence
between  solutions of the Nehari-Takagi problem and solutions of
some Takagi-Sarason interpolation problem. In Section 6 we calculate
the resolvent matrix of the Nehari-Takagi problem in the rational
case and show that it coincides with the resolvent matrix found
in~\cite{BGR90}.

\section{Preliminaries }\addcontentsline{toc}{section}{Preliminaries }
\subsection{Notations}
  Let  $\Omega_+$ be either   $\mathbb{D}=\{\lambda\in\mathbb{C}:\,|\lambda|<1\}$ or $\mathbb{C}_+=\{\lambda\in\mathbb{C}:\,\Im\lambda|>0\}$. Let us set for arbitrary $\lambda,\omega\in\dC$
\[
\rho_{\omega}(\lambda)=\left\{\begin{array}{cc}
                                1-\lambda\bar\omega, & \Omega_+=\dD, \\
                                - i(\lambda-\bar\omega), & \Omega_+=\dC_+,
                              \end{array}
\right.\quad
\lambda^\circ=\left\{\begin{array}{cc}
                                1/{\bar\lambda}, & \Omega_+=\dD, \\
                                \bar\lambda, & \Omega_+=\dC_+.
                              \end{array}
\right.
\]
Thus, $\Omega_+=\{\omega\in \mathbb{C}:\rho_\omega(\omega)>0\}$ and let
\[
\Omega_0=\{\omega\in \mathbb{C}:\rho_\omega(\omega)=0\},\quad \Omega_- =
\{\omega\in\mathbb{C }: \rho_{\omega}(\omega) < 0\}.
\]
The following basic classes of mvf's will be used in this paper:
 $H_2^{\ptq}$ (resp., $H_\infty^{\ptq}$) is  the
class of $\ptq$ mvf's with entries in the Hardy space $H_2$ (resp.,
$H_\infty$); $H_2^p:=H_2^{p\times 1}$,   and
$(H_2^p)^\perp=L_2^p\ominus H_2^p$,
$\cS^{\ptq}$ is the Schur class of $\ptq$ mvf's holomorphic and contractive on $\Omega_+$,
$\mathcal{S}_{in}^{p\times q}$ (resp., $\mathcal{S}_{out}^{p\times q}$) is the class of inner (resp., outer) mvf's in $\mathcal{S}^{p\times q}$:
\[
 \begin{split}
 &\mathcal{S}_{in}^{p\times q}=\{s\in \mathcal{S}^{p\times q}: s(\mu)^*s(\mu)=I_p \ \mbox{a.e. on}\ \Omega_0 \};\\
&\mathcal{S}_{out}^{p\times q}=\{s\in \mathcal{S}^{p\times q}:
\overline{s H_2^q}=H_2^p
\},
 \end{split}
\]
The Nevanlinna class $\mathcal{N}^{p\times q}$ and  the Smirnov class $\mathcal{N}^{p\times q}_{+}$ are defined by
\begin{equation}\label{eq:Smirn}
 \begin{split}
&\mathcal{N}^{p\times q}=\{f=h^{-1}g: \,g\in H_{\infty}^{p\times q},\, h\in  \mathcal{S}:=\mathcal{S}^{1\times 1}\},\\
&\mathcal{N}^{p\times q}_{+}=\{f=h^{-1}g: \,g\in H_{\infty}^{p\times q},\, h\in \mathcal{S}_{out}:=\mathcal{S}_{out}^{1\times 1}\}
.\\
 \end{split}
\end{equation}

For a mvf $f(\lambda)$ let us set
$
f^{\#}(\lambda)=f(\lambda^\circ)^*.
$
Denote by $\mathfrak{h}_f$ the domain of holomorphy of the mvf $f$ and let
$ \mathfrak{h}_f^\pm=\mathfrak{h}_f\cap\Omega_\pm$.

A $p\times q$ mvf $f_-$ in $\Omega_-$ is said to be a pseudocontinuation of a mvf $f\in\mathcal{N}^{p\times q}$, if
\begin{enumerate}
 \item [(1)] $f_{-}^\#\in\mathcal{N}^{p\times q}$;
     \item [(2)] $\lim_{\nu\downarrow 0}f_-(\mu-i\nu)=\lim_{\nu\downarrow 0}f_+(\mu+i\nu)(=f(\mu))$ \ a.e. on $\Omega_0.$

\end{enumerate}
The subclass of all mvf's $f\in\mathcal{N}^{p\times q}$ that admit pseudocontinuations
$f_-$ into $\Omega_-$ will be denoted $\Pi^{p\times q}$.

Let  $\varphi(\lambda)$ be a $p\times q$ mvf that is meromorphic on $\Omega_+$ with
a Laurent expansion
\[
    \varphi(\lambda)=(\lambda-\lambda_0)^{-k}\varphi_{-k}+\cdots
+(\lambda-\lambda_0)^{-1}\varphi_{-1}+\varphi_0+\cdots
\]
 in a neighborhood of its pole $\lambda_0\in\Omega_+$. The pole
multiplicity $\cM_{\pi}(\varphi, \lambda_0)$ is defined by (see~\cite{KL})
\[
\cM_{\pi}(\varphi, \lambda_0)=\rank L(\varphi,\lambda_0), \quad
T(\varphi,\lambda_0)=\begin{bmatrix}
  \varphi_{-k} &  & {\bf 0} \\
  \vdots & \ddots &  \\
  \varphi_{-1} &\hdots & \varphi_{-k}
\end{bmatrix}.
\]
The pole multiplicity of $\varphi$ over $\Omega_+$ is given by
\[
  \cM_{\pi}(\varphi, \Omega_+)=\sum_{\lambda\in\Omega_+}\cM_{\pi}(\varphi, \lambda).
\]
This definition of pole multiplicity coincides with that based on the
Smith-McMillan representation of $\varphi$ (see \cite{BGR90}).

Let $b_\omega(\lambda)$,  be a Blaschke factor
($b_\omega(\lambda)=\frac{\lambda-\omega}{1-\lambda\omega^*}$ in the case
$\Omega_+=\mathbb{D}$, $b_\omega(\lambda)=\frac{\lambda-\omega}{\lambda-\omega^*}$ in
the case $\Omega_+=\mathbb{C}_+$), and let $P$ be an orthogonal projection in $\mathbb{C}^p$.
Then the mvf
\[
B_\alpha (\lambda)=I_p-P+b_\alpha(\lambda)P, \quad\omega\in\Omega_+,
\]
belongs to the Schur class $\mathcal {S}^{p\times p}$ and is called \textit{the
elementary Blaschke--Potapov (BP) factor} and $B(\lambda)$ is called \textit{primary} if
$rank \ P=1$.
 The product
\[
B(\lambda)=\prod
\limits_{j=1}^{\overset{\kappa}{\curvearrowright}}B_{\alpha_j}(\lambda),
\]
 where
$B_{\alpha_j}(\lambda)$ are primary Blaschke--Potapov factors is called \textit{a
Blaschke--Potapov product} of degree $\kappa$.

\begin{remark}\label{rem:3.1}
For a Blaschke-Potapov product $b$
 the following statements are equivalent:
\begin{enumerate}
  \item[(1)] the degree of  $b$ is equal $\kappa$;
  \vskip 6pt
  \item[(2)] $\cM_{\pi}(b^{-1}, \Omega_+)=\kappa$;\vskip 6pt
\end{enumerate}
\end{remark}

\subsection{The generalized Schur class}

Let $\kappa \in\mathbb{Z}_{+}:=\mathbb{N}\cup\{0\}$. Recall, that a Hermitian
kernel
$\mathsf{K}_{\omega}(\lambda):\Omega\times\Omega\to\mathbb{C}^{m\times
m}$ is said to have $\kappa$ negative squares, if for every positive
integer $n$ and every choice of $\omega_j\in\Omega$ and
$u_j\in\mathbb{{C}}^m$ $(j=1,\dots,n)$ the matrix
$$({\left<{\mathsf K}_{\omega_j}(\omega_k)u_j,u_k\right>})_{j,k=1}^n$$
 has at most $\kappa$
negative eigenvalues, and for some choice of
$\omega_1,\ldots,\omega_n\in\Omega$ and $u_1,\ldots,u_n\in\mathbb{C}^m$
exactly $\kappa$ negative eigenvalues  (see~\cite{KL}).

Let ${\mathcal S}_{\kappa}^{q\times p}$ denote {\it the generalized
Schur class} of $q\times p$ mvf's $s$ that are meromorphic in
$\Omega_+$ and for which the kernel
\begin{equation}\label{kerLambda}
{\mathsf \Lambda}_\omega^s(\lambda)=
\frac{I_{p}-s(\lambda)s(\omega)^*}{\rho_\omega(\lambda)}
\end{equation}
has $\kappa$ negative squares on ${\mathfrak h}_s^+\times{\mathfrak h}_s^+$.
In the case where $\kappa=0$, the class ${\mathcal S}_{0}^{q\times
p}$ coincides with the Schur class ${\mathcal S}^{q\times p}$ of
contractive
 mvf's holomorphic in $\Omega_+$.
As was shown in~\cite{KL} every mvf $s\in {\mathcal S}_{\kappa}^{q
\times p}$ admits factorizations of the form
\begin{equation}\label{KLleft}
s(\lambda)=b_{\ell}(\lambda)^{-1}s_{\ell}(\lambda)=s_r(\lambda)b_r(\lambda)^{-1}, \quad  \lambda\in\mathfrak{h}_s^+,
\end{equation}
where $b_{\ell}\in {\mathcal S}^{q\times q}$,  $b_r\in {\mathcal S}^{p \times p}$ are  Blaschke--Potapov products of degree $\kappa$, $s_{\ell}, s_r\in{\mathcal S}^{q \times p}$
and the factorizations~\eqref{KLleft} are left coprime  and right coprime, respectively, i.e.
\begin{equation}\label{KLcanon2}
\textmd{rank} \left[\begin{array}{cc}
 b_{\ell}(\lambda) & s_{\ell}(\lambda)
\end{array}\right]
=q\quad (\lambda\in\Omega_+)
\end{equation}
and
\begin{equation}\label{KRcanon2.2}
\textmd{rank} \left[\begin{array}{cc}
b_r(\lambda)^* & s_r(\lambda)^*
\end{array}\right]  =p\quad (\lambda\in\Omega_+).
\end{equation}

The following matrix identity was established in the rational case in~\cite{Fran87}, in general case see~\cite{DD09}.
 \begin{theorem}\label{2.1}
\label{thm:jun1a8} Let $s \in{\mathcal S}^{q\times p}_\kappa$ have Kre\u{\i}n-Langer
factorizations
\begin{equation}\label{KL_1}
s=b_\ell^{-1}s_\ell =s_rb_r^{-1}.
\end{equation}
Then there exists a set of mvf's  $c_\ell\in H_\infty^{q\times q}$, $d_\ell\in
H_\infty^{p \times q}$, $c_r\in H_\infty^{p \times p}$ and $d_r\in H_\infty^{p\times q}$,
such that
\begin{equation}
\label{eq:jun1a8}
\begin{bmatrix}c_r&d_r\\-s_\ell&b_\ell\end{bmatrix}
\begin{bmatrix}b_r&-d_\ell\\ s_r&c_\ell\end{bmatrix}=\begin{bmatrix}I_p&0\\
0&I_q\end{bmatrix}.
\end{equation}
\end{theorem}
\subsection{The generalized Smirnov class}
Let $\mathcal{R}^{p\times q}$ denote the class of rational $p\times q$ mvf's and let $\kappa\in\mathbb{N}$.
A $p\times q$ mvf $\varphi(z)$ is said to belong to the generalized Smirnov class $\mathcal{N}_{+,\kappa}^{p\times q}$, if it admits the representation
\[
\varphi(z)=\varphi_0(z)+r(z), \quad\mbox{where}\quad
\varphi_0\in\mathcal{N}_+^{p\times q}, \, r\in\mathcal{R}^{p\times q} \quad \text{and}\quad M_\pi(r,\Omega_+)\leq\kappa.
\]
If $\kappa=0$,  then the class $\mathcal{N}_{+,0}^{p\times q}$ coincides with the Smirnov class $\mathcal{N}_{+}^{p\times q}$, defined in~\eqref{eq:Smirn}.
The generalized Smirnov class $\mathcal{N}_{+,\kappa}^{p\times q}$ was introduced in~\cite{Nei08}. In~\cite{DD14},  mvf's $\varphi$ from $\mathcal{N}_{+,\kappa}^{p\times q}$ were characterized by  the following left coprime factorization
 \[
    \varphi(\lambda)=b_\ell(\lambda)^{-1}\varphi_\ell(\lambda),
 \]
where $b_\ell\in S_{in}^{\ptp}$ is a Blaschke--Potapov product of degree $\kappa$,  $\varphi_\ell\in \cN_+^{\ptq}$ and
\[
    \rank \left[\begin{array}{cc}
b_{\ell}(\lambda) &\varphi_{\ell}(\lambda)
\end{array}\right]
=p\quad  \textrm{for}\ \lambda\in\Omega_+.
\]
Clearly, for $\varphi\in\mathcal{N}_{+,\kappa}^{p\times q}$ there exists a right coprime factorization with similar  properties. This implies, in particular, that the class $\cS_\kappa^{p\times q}$ is contained in $\mathcal{N}_{+,\kappa}^{p\times q}$.
\subsection{Generalized $j_{pq}$-inner mvf's}
Let $j_{pq}$  be an $m\times m$ signature matrix
$$j_{pq}=\left[\begin{array}{cc}
      I_p & 0 \\
      0 &-I_q   \end{array} \right],\quad\text{where}\ p+q=m,$$

\begin{definition}\label{U_kappa}~\cite{AlpD86}
An $m\times m$ mvf  $W(\lambda)=[w_{ij}(\lambda)]_{i,j=1}^2$ that is
meromorphic in $\Omega_+$ is said to belong to the class ${\mathcal
U}_{\kappa}(j_{pq})$ of {\it generalized $j_{pq}$-inner} mvf's, if:
\begin{enumerate}
\item[(i)]
the kernel
\[
{\mathsf K}_\omega^W(\lambda)= \frac{j_{pq}-W(\lambda)j_{pq}W(\omega)^*}{\rho_\omega(\lambda)}
\]
has $\kappa$ negative squares in ${\mathfrak h}_W^+\times{\mathfrak h}_W^+$;
\item[(ii)]
$j_{pq}-W(\mu)j_{pq}W(\mu)^*=0$ a.e. on $\Omega_0$.
\end{enumerate}
\end{definition}

As is known \cite[Th.6.8.]{AlpD86} for every $W\in\mathcal{U}_\kappa(j_{pq})$ the block $w_{22}(\lambda)$ is invertible for all $\lambda\in\mathfrak{h}^+_W$ except for at most $\kappa$  points in $\Omega_+$. Thus the Potapov-Ginzburg  transform of $W$
\begin{equation}\label{PGtrans2}
    S(\lambda)=PG(W):=\left[\begin{array}{cc}
      w_{11}(\lambda) & w_{12}(\lambda) \\
      0 & I_q
    \end{array}      \right]
    \left[\begin{array}{cc}
            I_p &       0\\
      w_{21}(\lambda) & w_{22}(\lambda)
    \end{array}      \right]^{-1}
\end{equation}
is well defined for those $\lambda\in{\mathfrak h}_W^+$, for which
$w_{22}(\lambda)$ is invertible. It is well known that $S(\lambda)$
belongs to the class ${\mathcal S}_\kappa^{m\times m}$ and $S(\mu)$
is unitary for a.e. $\mu\in \Omega_0$ (see~\cite{AlpD86},
~\cite{DD09}).

\begin{definition}\label{eq:11.100}~\cite{DD09}
An $m\times m$ mvf $W\in\mathcal{U}_\kappa(j_{pq})$ is said to be in the class
$\mathcal{U}_\kappa^r(j_{pq})$, if
\begin{equation}\label{s21}
s_{21}:=-w_{22}^{-1}w_{21} \in\  {\mathcal S}_\kappa^{q\times p}.
\end{equation}
\end{definition}
Let $W\in\mathcal{U}_\kappa^r(j_{pq})$ and let the Krein-Langer factorization
of $s_{21}$ be written as
\[
    s_{21}(\lambda)={b}_{\ell}(\lambda)^{-1}{s}_{\ell}(\lambda)
={s}_r(\lambda){b}_r(\lambda)^{-1}\quad   (\lambda\in\mathfrak
{h}_{s_{21}}^+),
\]
where ${b}_{\ell}\in {\mathcal S}^{q\times q}_{in}$, $b_{r}\in {\mathcal S}^{p\times
p}_{in}$, ${s}_{\ell},{s}_r\in {\mathcal S}^{q\times p}$. Then, as was shown
in~\cite{DD09}, the mvf's  $b_{\ell}s_{22}$ and $s_{11}{b}_{r}$ are holomorphic in
$\Omega_+$, and
\[
{b}_{\ell}s_{22}\in \mathcal{S}^{q\times q} \quad \text{and} \quad s_{11}{b}_r\in
\mathcal{S}^{p\times p}.
\]
\begin{definition}\label{eq:ap}~\cite{DD09}
Consider inner-outer factorization of $s_{11}{b}_r$ and outer-inner factorization of $b_{\ell}s_{22}$
\begin{equation}\label{eq:15}
s_{11}{b}_r=b_1 a_1, \qquad {b}_{\ell}s_{22}=a_2 b_2,
\end{equation}
where $b_1\in {\mathcal S}_{in}^{p\times p}$, $b_2\in {\mathcal S}_{in}^{q\times q}$,
$a_1\in {\mathcal S}_{out}^{p\times p}$, $a_2\in {\mathcal S}_{out}^{q\times
q}$. The pair $\{b_1,b_2\}$ of inner factors in the factorizations~\eqref{eq:15} is called
{\it the associated pair} of the mvf $W\in\mathcal{U}_\kappa^r(j_{pq})$.

\end{definition}
From now onwards this pair $\{b_1,b_2\}$ will be called also a right associated pair since it is related to the right linear fractional transformation
\begin{equation}\label{eq:lft}
    T_W[\varepsilon]:=(w_{11}\varepsilon+
w_{12})(w_{21}\varepsilon+w_{22})^{-1},
\end{equation}
see~\cite{Arov88},~\cite{ArovD97,ArovD08}. Such transformations play
important role in description of solutions of different interpolation
problems, see~\cite{AAK71a}, \cite{Arov88}, \cite{BGR90}, \cite{AI},
 \cite{Der03}, \cite{DD10}. In the case $\kappa=0$ the
definition of the associated pair was given in~\cite{Arov88}.

For every $W\in\mathcal{U}_\kappa^r(j_{pq})$ and $\varepsilon\in\cS^{\ptq}$ the mvf
$T_W[\varepsilon]$ admits the dual  representation
$$T_W[\varepsilon]=(w_{11}^\#+\varepsilon
w_{12}^\#)^{-1}(w_{21}^\#+\varepsilon w_{22}^\#).$$
As was shown in~\cite{DD09}, for $W\in\mathcal{U}_\kappa^r(j_{pq})$ and $c_r$, $d_r$, $c_\ell$ and $d_\ell$ as in~\eqref{eq:jun1a8} the mvf
\begin{equation}\label{eq:jun10a9}
K^\circ:=(-w_{11}d_\ell+w_{12}c_\ell)(-w_{21}d_\ell+w_{22}c_\ell)^{-1},
\end{equation}
 belongs to $H_\infty^{p\times q}$.
It is clear that
 $(K^\circ)^\# \in H^{q\times p}_\infty (\Omega_-)$.

 In the future we will need the following factorization of the mvf $W\in
\mathcal{U}_\kappa^r(j_{pq})$, which was obtained
in~\cite[Theorem~4.12]{DD09}:
\begin{equation}
\label{eq:11.34}
W=\Theta^\circ\,\Phi^\circ\quad in\  \Omega_+,
\end{equation}
where
\[
    \Theta^\circ=\begin{bmatrix}b_1&K^\circ b_2^{-1}\\0&b_2^{-1}\end{bmatrix}, 
    \quad
    \Phi^\circ,\, (\Phi^\circ)^{-1} \in \cN_+ .
\]

\section{The Takagi-Sarason interpolation problem}
\textbf{Problem $\mathbf{TSP}_\kappa(b_1,b_2,K)$ }
Let
$b_1\in \cS_{in}^{p\times p}$, $b_2\in \cS_{in}^{q\times q}$ be
inner mvf's,
let $K\in H_\infty^{p\times q}$ and let $\kappa\in \dZ$.
A ${p\times q}$ mvf $s$ is called a solution of the Takagi-Sarason
problem $\mathbf{TSP}_\kappa(b_1,b_2,K)$, if  $s$  belongs to
$\cS_{\kappa'}^{p\times q}$ for some $\kappa'\le \kappa$ and
satisfies~
\begin{equation}\label{eq:3.15}
b_1^{-1}(s-K)b_2^{-1}\in\mathcal{N}_{+,\kappa}^{p\times q}.
\end{equation}

The set of solutions of the Takagi-Sarason
problem will be denoted by
$$\mathcal{TS}_\kappa(b_1, b_2, K)=\bigcup_{\kappa'\leq \kappa}\{s\in\mathcal{S}_{\kappa'}^{p\times q}: b_1^{-1}(s-K)b_2^{-1}\in\mathcal{N}_{+, \kappa}^{p\times q}\}.$$
The problem $\mathbf{TSP}_\kappa(b_1,b_2,K)$ has been studied in~\cite{BH83}, in the rational case ($K\in\cR^{q\times q}$) the set  $\mathcal{TS}_\kappa(b_1, b_2, K)\cap\cR^{p\times q}$ was described in~\cite{BGR90}. In the completely indeterminate case explicit formulas for the resolvent matrix  can be found in~\cite{DD10}, \cite{DD14}. In the general positive semidefinite case, the problem was solved in~\cite{KKhYu87}, \cite{Kh90}.

We now recall the construction of the resolvent matrix  from \cite{DD14}.
Let
$$
\cH(b_1)=H_2^p\ominus b_1H_2^p,
\quad
\cH_*(b_2):=(H_2^q)^\perp\ominus
b_2^*(H_2^q)^\perp
$$
let
$$
\cH(b_1,b_2):=\cH(b_1)\oplus\cH_*(b_2).
$$
and let the operators $ K_{11}: H_2^q\to\cH(b_1)$, $K_{12}:
\cH_*(b_2)\to\cH(b_1)$, $K_{22}: \,\cH_*(b_2)\to(H_2^p)^\perp $ and $P: \cH(b_1,b_2)\to\cH(b_1,b_2)$
be defined by the formulas
\begin{equation}\label{eq:I26}
\begin{split}
K_{11}h_+&=\Pi_{\cH(b_1)}Kh_+,\quad h_+\in H_2^q,\\
K_{12}h_2&=\Pi_{\cH(b_1)}Kh_2,\,\quad h_2\in\cH_*(b_2),\\
K_{22}h_2&=\Pi_-Kh_2,\,\,\,\qquad h_2\in\cH_*(b_2),
\end{split}
\end{equation}
\begin{equation}\label{eq:I30}
    P=\begin{bmatrix}
  I-K_{11}K_{11}^* & -K_{12} \\
  -K_{12}^* & I -K_{22}^*K_{22}\\
\end{bmatrix}.
\end{equation}
The data set $b_1,\,b_2,\,K$ considered in \cite{DD14} is subject to the following constraints:
\begin{enumerate}
    \item [(H1)] $b_1\in \cS^{\ptp}_{in}$, $b_2\in
    \cS^{\qtq}_{in}$, \,$K\in H^{\ptq}_\infty$.
\vspace{2mm}
    \item [(H2)] $\kappa_1=\nu_-(P)<\infty$.
\vspace{2mm}
    \item [(H3)]
$0\in\rho(P)$.
\vspace{2mm}
    \item [(H4)]  $\gh_{b_1}\cap\gh_{b_2^\#}\cap\Omega_0\ne \emptyset$.
 \end{enumerate}

Define also the operator
\begin{equation}\label{eq:I31}
    {F}=\begin{bmatrix}
  I & K_{22} \\
  K_{11}^* & I \\
\end{bmatrix}\,:\begin{array}{c}\cH(b_1)\\ \oplus\\ \cH_*(b_2)\end{array}
\rightarrow \begin{array}{c}b_1(H_2^p)^\perp\\ \oplus\\ b_2^*(H_2^q)\end{array}\stackrel{def}{=}\cK.
\end{equation}
As was shown in \cite{DD14} for every $h_1\in\cH(b_1)$ and $h_2\in\cH_*(b_2)$ the vvf's $(K_{11}^*h_1)(\lambda)$ and $(K_{22}h_2)(\lambda)$
  admit pseudocontinuations of bounded type which are holomorphic on
  ${\gh_{b_1}}$ and ${\gh_{b_2^\#}}$, respectively.
This allows to define an $m\times m$ mvf  $\lambda\to F(\lambda)$  by  
\begin{equation}
\label{eq:oct4a11}
F(\lambda)=E(\lambda)F\quad\textrm{for }\lambda\in
\gh_{b_1}\cap\gh_{b_2^\#}
\end{equation}
where $E(\lambda)$ is the evaluation operator
$$
E(\lambda)\,:f\in\cK\to f(\lambda)\in\dC^m.
$$
Let $\mu\in\gh_{b_1}\cap\gh_{b_2^\#}\cap\Omega_0.
$ Then the mvf $W(\lambda)$ defined by
\begin{equation}\label{eq:2.11}
    W(\lambda)=I-\rho_\mu(\lambda)F(\lambda)P^{-1}F(\mu)^*j_{pq}\quad
\textrm{for $\lambda\in\gh_{b_1}\cap\gh_{b_2^\#}$}
\end{equation}
belongs to the class $\cU_{\kappa_1}^r(j_{pq})$ of {\it generalized $j_{pq}$-inner} mvf's and takes values in $L_2^{m\times m}$.
The following theorem presents a description of the set $\mathcal{TS}_\kappa(b_1, b_2, K)$.
\begin{theorem}
\label{thm:nov11a13}
Let (H1)--(H4) be in force and let $W(\lambda)$ be the mvf, defined by~\eqref{eq:2.11}. Then $W\in\cU_{\kappa_1}^r(j_{pq})\cap L_2^{m\times m}$ and
\begin{enumerate}
\item[\rm(1)]  $\cT\cS_\kappa(b_1,b_2;\,K)\ne\emptyset \Longleftrightarrow \nu_-(P)\le \kappa$.
\vspace{2mm}
\item[\rm(2)] If  $\kappa_1=\nu_-(P)\le \kappa$, then
\begin{equation}\label{eq:ST_descr}
\cT\cS_\kappa(b_1,b_2;K)
=T_W[{\cS}_{\kappa-\kappa_1}^{\ptq}]
:=\{T_W[\varepsilon]:\,\varepsilon\in\cS_{\kappa-\kappa_1}^{\ptq}\},
\end{equation}
where $T_W[\varepsilon]$ is the linear fractional transformation given by~\eqref{eq:lft}.
\end{enumerate}
\end{theorem}
\begin{proof}
The proof of this statement can be derived from the proof of Theorem~5.7 in~\cite{DD14}.
However, we would like to present here a shorter proof based on the description of the set $\cT\cS_\kappa(b_1,b_2;\,K)$,
given in~\cite[Theorem~5.17]{DD10}.

As was shown in \cite[see Theorem~4.2 and Corollary~4.4]{DD14}
the mvf $W(z)$ belongs to the class $\cU_{\kappa_1}^r(j_{pq})$ of  generalized $j_{pq}$-inner mvf's
with the property~\eqref{s21} and $\{b_1,b_2\}$ is the associated pair of $W$.
Moreover, by construction  $W(z)$ takes values in $L_2^{m\times m}$.
Let $c_\ell$ and $d_\ell$ be mvf's defined in Theorem~\ref{thm:jun1a8} and let  $K^{\circ}$ be given
by~\eqref{eq:jun10a9}. Then $W$ admits the factorization~\eqref{eq:11.34} (see~\cite[Theorem~4.12]{DD09}).
This proves that all the assumptions of Theorem~5.17  from~\cite{DD10}
with $K$ replaced by $K^{\circ}$ are satisfied and by that theorem
\begin{equation}\label{eq:STdescr}
    {\cT\cS}_\kappa(b_1,b_2;K^{\circ})
=T_W[{\cS}_{\kappa-\kappa_1}^{\ptq}].
\end{equation}

On the other hand it follows from \cite[Theorem~4.2]{DD14} that the mvf $W$  admits the factorization
\begin{equation}\label{eq:W_fact}
    W=\Theta\,\Phi=\begin{bmatrix}b_1&Kb_2^{-1}\\0&b_2^{-1}\end{bmatrix}
\begin{bmatrix}\varphi_{11}& \varphi_{12}\\
\varphi_{21}&\varphi_{22}\end{bmatrix}
\end{equation}
with $\Phi,\Phi^{-1}\in\cN_+^{m\times m}$.
Comparing~\eqref{eq:W_fact} with~\eqref{eq:11.34} one obtains
\[
\begin{bmatrix}I&b_1^{-1}(K-K^\circ)b_2^{-1}\\0&I\end{bmatrix}
=\Phi^\circ\Phi^{-1}\in\cN_+^{m\times m}
\]
and hence
\begin{equation}\label{eq:KK0}
    b_1^{-1}(K-K^\circ)b_2^{-1}\in\cN_+^{p\times q}.
\end{equation}
This implies the equality $\cT\cS_\kappa(b_1,b_2;K)=\cT\cS_\kappa(b_1,b_2;K^\circ)$,
that in combination with~\eqref{eq:STdescr} completes the proof.
\end{proof}
\begin{remark}\label{rem:GSTP}
Alongside with the set $\mathcal{TS}_\kappa(b_1, b_2, K)$ consider also its subset
\begin{equation}\label{eq:Sdef}
\mathcal{S}_\kappa(b_1, b_2, K)=\{s\in\mathcal{S}_{\kappa}^{p\times q}: b_1^{-1}(s-K)b_2^{-1}\in\mathcal{N}_{+, \kappa}^{p\times q}\}.
\end{equation}
A ${p\times q}$ mvf $s$ is called a solution of the  generalized Schur-Takagi problem \linebreak $\mathbf{GSTP}_\kappa(b_1,b_2,K)$ if  $s$  belongs to the set $\mathcal{S}_\kappa(b_1, b_2, K)$. A description of the set $\mathcal{S}_\kappa(b_1, b_2, K)$ was obtained in~\cite[Theorem~1.2]{DD14} in the form:
\begin{equation}\label{eq:GSTP_kappa}
    \mathcal{S}_\kappa(b_1, b_2, K)=
    T_W[{\cS}_{\kappa-\kappa_1}^{\ptq}]\cap S_\kappa^{\ptq}.
\end{equation}
Notice, that the reasonings of Theorem~\ref{thm:nov11a13} allows to give a shorter proof of the formula~\eqref{eq:GSTP_kappa}. Indeed, by~\cite[Theorem~5.17]{DD10}
\begin{equation}\label{eq:Sdescr}
    {\cS}_\kappa(b_1,b_2;K^{\circ})
=T_W[{\cS}_{\kappa-\kappa_1}^{\ptq}]\cap \cS_\kappa^{\ptq}.
\end{equation}
Next, it follows from~\eqref{eq:KK0} that
\begin{equation}\label{eq:GS_kappa}
\cS_\kappa(b_1,b_2;K)=\cS_\kappa(b_1,b_2;K^\circ).
\end{equation}
Now~\eqref{eq:GSTP_kappa} is implied by~\eqref{eq:GS_kappa} and~\eqref{eq:Sdescr}.
\end{remark}
\section{Generalized $\gamma$-generating mvf's}
\begin{definition}\label{gamma_r}
Let $\mathfrak{M}^r_\kappa(j_{pq})$ denote the class of $m\times m$ mvf's $\mathfrak{A}(\mu)$ on $\Omega_0$ of the form
 $$\mathfrak{A}(\mu)=\left[\begin{array}{cc}
      a_{11}(\mu) & a_{12}(\mu) \\
       a_{21}(\mu) &  a_{22}(\mu)
    \end{array}      \right],$$
    with blocks $a_{11}$ and $a_{22}$ of size $p\times p$ and $q\times q$, respectively,  such that:
    \begin{enumerate}
\item [(1)]  $\mathfrak{A}(\mu)$ is a measurable mvf on $\Omega_0$ and $j_{pq}$-unitary a.e. on $\Omega_0$;
 \item[(2)] the mvf's $a_{22}(\mu)$ and $a_{11}(\mu)^{*}$ are invertible for a.e. $\mu\in\Omega_0$ and the mvf
 \begin{equation}\label{eq:s_21}
    s_{21}(\mu)=-a_{22}(\mu)^{-1}a_{21}(\mu)
                =-a_{12}(\mu)^{*}(a_{11}(\mu)^{*})^{-1}
\end{equation}
     is the boundary value of a mvf $s_{21}(\lambda)$ that belongs to the class $\mathcal{S}_\kappa^{q\times p}$;
 \item[(3)] $a_{11}(\mu)^{*}$ and  $a_{22}(\mu)$, are the boundary values of mvf's $a_{11}^\#(\lambda)$ and
$a_{22}(\lambda)$ that are meromorphic in $\dC_+$ and, in addition,
\begin{equation}\label{eq:a_1a_2}
    a_1:=(a_{11}^\#)^{-1}b_r\in\mathcal{S}_{out}^{p\times p},\quad
    a_2:= b_\ell a_{22}^{-1}\in\mathcal{S}_{out}^{q\times q},
\end{equation}
  where
$b_\ell$, $b_r$ are Blaschke-Potapov products of degree $\kappa$, determined by  Krein-Langer factorizations of
$s_{21}$.
\end{enumerate}
\end{definition}
Mvf's in the class $\mathfrak{M}^r_\kappa(j_{pq})$ are called generalized right $\gamma$-generating mvf's.
The class $\mathfrak{M}^r(j_{pq}):=\mathfrak{M}^r_0(j_{pq})$ was introduced in~\cite{Arov89},
in this case conditions (2) and (3) in Definition~\ref{gamma_r} are simplified to:
    \begin{enumerate}
 \item[(2$'$)] $s_{21}\in\mathcal{S}^{q\times p}$;
 \item[(3$'$)] $    a_1:=(a_{11}^\#)^{-1}\in\mathcal{S}_{out}^{p\times p}$,
    $a_2:=  a_{22}^{-1}\in\mathcal{S}_{out}^{q\times q}$.
\end{enumerate}
Mvf's from the class $\mathfrak{M}^r(j_{pq})$ play an important role in the description of solutions of the
Nehari problem and are called right $\gamma$-generating mvf's, \cite{Arov89,ArovD08}.
Mvf's in the class $\mathfrak{M}^r_\kappa(j_{pq})$ will be called generalized right $\gamma$-generating mvf's.
\begin{definition}~\cite{ArovD08}
An ordered pair $\{b_1, b_2\}$ of inner mvf's $b_1\in \mathcal{S}^{p\times p}$, $b_2\in \mathcal{S}^{q\times q}$
is called a denominator of the mvf $f\in\mathcal{N}^{p\times q}$, if
\[
b_1fb_2\in\mathcal{N}^{p\times q}_+.
\]
The set of  denominators of the mvf $f\in\mathcal{N}^{p\times q}$ is denoted by $\mbox{den }(f)$.
\end{definition}

\begin{theorem}\label{thm:7.26}
  Let $\mathfrak{A}\in\Pi^{m\times m}\cap\mathfrak{M}_\kappa^r(j_{pq})$, and let $c_r$, $d_r$, $c_\ell$ and $d_\ell$ be as in
Theorem~\ref{thm:jun1a8}.
 \begin{equation}\label{eq:4.38}
 f_0=(-a_{11}d_\ell+a_{12}c_\ell)a_2.
 \end{equation}
  Then the mvf $f_0$ admits the dual representation
 \begin{equation}\label{eq:11.29}
    f_0=a_1(c_ra_{21}^\#-d_ra_{22}^\#).
\end{equation}
If, in addition,
 $\{b_1,b_2\}\in den(f_0)$  and
  \begin{equation}\label{eq:4.39}
W(z)=\begin{bmatrix} b_1&0\\0&b_2^{-1}   \end{bmatrix}\mathfrak{A}(z),
 \end{equation}
then
$W\in\mathcal{U}_\kappa^r(j_{pq})$ and $\{b_1,b_2\}$ is the associated pair of $W$.

 Conversely, if $W\in\mathcal{U}_\kappa^r(j_{pq})$ and $\{b_1,b_2\}$ is the associated pair of $W$, then
 \[
     \mathfrak{A}(z)=
     \begin{bmatrix} b_1^{-1}&0\\0&b_2   \end{bmatrix}
     W(z)\in\Pi^{m\times m}\cap\mathfrak{M}_\kappa^r(j_{pq}) \quad \mbox{and}\quad \{b_1,b_2\}\in den(f_0).
 \]
 \end{theorem}

 \begin{proof}
 Let $\mathfrak{A}\in\Pi^{m\times m}\cap\mathfrak{M}_\kappa^r(j_{pq})$.
 It follows from~\eqref{eq:s_21}, \eqref{eq:a_1a_2} and~\eqref{KLleft} that
\[
\begin{split}
    -a_{21}d_\ell+a_{22}c_\ell
       &=\begin{bmatrix}
              a_{21}& a_{22}
       \end{bmatrix}
       \begin{bmatrix}
              -d_\ell \\
             c_\ell
       \end{bmatrix} =
       \begin{bmatrix}
              -a_{22}s_{21}& a_{22}
       \end{bmatrix}
       \begin{bmatrix}
              -d_\ell \\
             c_\ell
       \end{bmatrix} \\
       &=a_{22}b_{\ell}^{-1}
       \begin{bmatrix}
              -s_\ell& b_\ell
       \end{bmatrix}
       \begin{bmatrix}
              -d_\ell \\
             c_\ell
       \end{bmatrix}\\
       &=a_2^{-1}(s_\ell d_\ell+b_\ell c_\ell)=a_2^{-1}.
\end{split}
\]
 Let $f_0$ be defined by the equation~\eqref{eq:4.38}.
Then
 \[
     f_0=(-a_{11}d_\ell+a_{12}c_\ell)(-a_{21}d_\ell+a_{22}c_\ell)^{-1}.
 \]
The identity
\[
    \begin{bmatrix}
      c_r & -d_r
    \end{bmatrix}
     \mathfrak{A}^\#j_{pq} \mathfrak{A}
    \begin{bmatrix}
      -d_\ell\\ c_\ell
    \end{bmatrix}
    =\begin{bmatrix}
      c_r & -d_r
    \end{bmatrix}
    j_{pq}
    \begin{bmatrix}
      -d_\ell\\ c_\ell
    \end{bmatrix}=0
\]
implies that
\[
    (c_r a_{11}^\#-d_ra_{12}^\#)(-a_{11}d_\ell+a_{12}c_\ell) =
    (c_r a_{21}^\#- d_ra_{22}^\#)(-a_{21}d_\ell+a_{22}c_\ell),
\]
and hence that $f_0$ admits the dual representation
\[
    f_0=(c_ra_{11}^\#-d_ra_{12}^\#)^{-1}(c_r
    a_{21}^\#-d_ra_{22}^\#).
\]
Using the identity
\[
    \begin{split}
        \begin{bmatrix}
            c_r & -d_r
        \end{bmatrix}
         \begin{bmatrix}
              a_{11}^\#\\
              a_{12}^\#
        \end{bmatrix}&=
        \begin{bmatrix}
            c_r & -d_r
        \end{bmatrix}
         \begin{bmatrix}
              a_{11}^\#\\
              -s_{21}a_{11}^\#
        \end{bmatrix}\\
        &=\begin{bmatrix}
            c_r & -d_r
        \end{bmatrix}
         \begin{bmatrix}
              I_p\\
             -s_rb_{r}^{-1}
        \end{bmatrix}b_ra_1^{-1}=a_1^{-1}
    \end{split}
 \]
one obtains the equality
\[
    f_0=a_1(c_ra_{21}^\#-d_ra_{22}^\#)
\]
which coincides with~\eqref{eq:11.29}.

 Let $\{b_1, b_2\}\in den(f_0)$, i.e. $b_1f_0b_2\in \cN_+^{p\times q}$. Since $b_1f_0b_2\in L_\infty^{p\times q}$
 then by Smirnov theorem
 \[
    b_1f_0b_2\in H^{p\times q}_\infty.
 \]
Let us find the Potapov-Ginzburg transform $S=PG(W)$ of $W$, see~\eqref{PGtrans2}. The formula \eqref{eq:4.39}
implies that
\begin{eqnarray}\label{eq:4.42}
    s_{21}&=&-w_{22}^{-1}w_{21}=-a_{22}^{-1}a_{21}=-b_\ell^{-1}s_\ell,\\
    \label{eq:4.43}
    s_{22}&=&w_{22}^{-1}=a_{22}^{-1}b_2=b_\ell^{-1}a_2b_2,\\
    \label{eq:4.44}
    s_{11}&=&w_{11}^{-*}=b_1a_1a_1^{-1}b_1^{-1}w_{11}^{-*}\\
    &=&b_1a_1(c_ra_{11}^*-d_ra_{12}^*)b_1^{-1}w_{11}^{-*}\nonumber\\
    &=&b_1a_1(c_rw_{11}^*-d_rw_{12}^*)w_{11}^{-*}\nonumber\\
    &=&b_1a_1(c_r+d_rs_{21}),\nonumber\\
    s_{12}&=&-w_{11}^{-*}w_{21}^*=
    b_1a_1(c_rw_{11}^*-d_rw_{12}^*)w_{11}^{-*}w_{21}^*\label{eq:4.45}\\
    &=&b_1a_1(c_rw_{11}^*-d_rw_{22}^*+d_rs_{22})\nonumber\\
    &=&b_1f_0b_2+b_1a_1d_rs_{22}.\nonumber
\end{eqnarray}
The equalities \eqref{eq:4.42}-\eqref{eq:4.45} lead easily to the formula
 \begin{equation}\label{eq:4.46}
    \begin{split}
        S(z)&=\begin{bmatrix} b_1a_1c_r+b_1a_1d_rs_{21}& b_1f_0b_2+b_1a_1d_rs_{22}\\ s_{21}& s_{22}\end{bmatrix}\\&=\begin{bmatrix}b_1a_1c_r & b_1f_0b_2\\0&0  \end{bmatrix}
        +\begin{bmatrix}b_1a_1d_r\\I  \end{bmatrix}\begin{bmatrix}s_{21}& s_{22} \end{bmatrix}
        \\&=T(z)+\begin{bmatrix}b_1a_1d_r\\I  \end{bmatrix}b_\ell^{-1}\begin{bmatrix}-s_{\ell}& a_2b_2 \end{bmatrix},
    \end{split}
 \end{equation}
where $T(z)\in H_\infty^{m\times m}$. It follows from~\eqref{eq:4.46} that
\[
M_\pi(S, \Omega_+)\leq \kappa.
\]
On the other hand
\[
M_\pi(s_{21},\Omega_+)=M_\pi(-b_\ell^{-1}s_\ell,\Omega_+)=\kappa,
\]
and, consequently,
\[
    M_\pi(S, \Omega_+)= \kappa.
\]
Thus, $S\in\mathcal{S}_{\kappa}^{m\times m}$ and, hence, $W\in\mathcal{U}^r_\kappa(j_{pq})$.
\end{proof}

\section{A Nehari-Takagi problem}
Let $f\in L^{p\times q}_\infty$ and let $\Gamma(f)$ be the Hankel operator associated with $f_0$:

\begin{equation}\label{eq:3.1}
    \Gamma(f):=\Pi_-M_f|_{H_2^q},
\end{equation}
where $M_f$  denotes the operator of multiplication by $f$, acting from $L_2^q$ into $L_2^p$ and let $\Pi_-$ denote the orthogonal projection of $L_2^p$ onto $(H_2^p)^\perp$. The operator $\Gamma(f)$ is bounded as an operator from $H_2^q$ to $(H_2^p)^\perp$, moreover,
\[
    \|\Gamma(f)\|\leq\|f\|_{L^{p\times q}_\infty}.
\]

Consider the following Nehari-Takagi problem

\noindent
\textbf{Problem $\mathbf{NTP}_\kappa(f_0)$}: Given a mvf $f_0\in L_\infty^{p\times q}$.
Find $f\in L_\infty^{p\times q}$, such that
\begin{equation}\label{eq:3.10}
    \rank(\Gamma(f)-\Gamma(f_0))\leq\kappa\quad\mbox{and}\quad \|f\|_\infty\leq1.
\end{equation}
In the scalar case, the problem $\mathbf{NTP}_\kappa(f_0)$ has been
solved  by V.M.~Adamyan, D.Z.~Arov and M.G.~Kre\u{\i}n
in~\cite{AAK68} for the case $\kappa=0$  and in~\cite{AAK71a} for
arbitrary $\kappa\in\mathbb{N}$. In the matrix case a description of
solutions of the problem $\mathbf{NTP}_0(f_0)$  was obtained in the
completely indeterminate case by V.M.~Adamyan,~\cite{Ad73}, and in
the general positive-semidefinite case by
A.~Kheifets,~\cite{Kh2000}. The indefinite case
$(\kappa\in\mathbb{N})$ was treated in~\cite{BH83} (see
also~\cite{BGR90}, where an explicit formula for the resolvent
matrix was obtained in the rational case).

In what follows we confine ourselves  to the case when $\mbox{den }(f_0)\ne\emptyset$  and give a description of all solutions of the problem $\mathbf{NTP}_\kappa(f_0)$.
Let us set for arbitrary $f_0\in L_\infty^{p\times q}$
\[
    \mathcal{N}_{\kappa}(f_0)=\{f\in L_\infty^{p\times q}: f-f_0\in\mathcal{N}_{+,\kappa}^{p\times q}, \|f\|\leq 1 \}
\]
and let us denote the set of solutions of the problem $\mathbf{NTP}_\kappa(f_0)$ by
\[
    \mathcal{NT}_{\kappa}(f_0)=\{f\in L_\infty^{p\times q}: \rank(\Gamma(f)-\Gamma(f_0))\leq\kappa
    \,\,\mbox{ and }\,\, \|f\|\leq 1 \}.
\]
By Kronecker Theorem (\cite{Gant}), the condition $f-f_0\in\mathcal{N}_{+,\kappa}^{p\times q} $ is equivalent to
\[
\rank(\Gamma(f)-\Gamma(f_0))=\kappa,
\]
Therefore, the set $\mathcal{NT}_{\kappa}(f_0)$ 
is represented as
\begin{equation}\label{eq:3.13}
\mathcal{NT}_{\kappa}(f_0)=\bigcup_{\kappa'\leq \kappa}\mathcal{N}_{\kappa'}(f_0).
\end{equation}
In the following theorem relations between the set of solutions of the Nehari-Takagi problem and the set of
solutions of a Takagi-Sarason problem is established in the case when $\mbox{den }(f_0)\ne \emptyset $.
\begin{theorem}\label{thm:3.1}
Let $f_0\in L^{p\times q}_\infty$, $\Gamma=\Gamma(f_0)$, $\kappa\in\mathbb{Z}_+$, $\{b_1, b_2\}\in den(f_0)$ and
$K=b_1f_0b_2$. Then
\[
f\in\mathcal{N}_\kappa(f_0)\Leftrightarrow s=b_1fb_2\in \mathcal{TS}_\kappa(b_1, b_2, K_0).
\]

\end{theorem}
\begin{proof}
Let $f\in\mathcal{N}_\kappa(f_0)$. Then the mvf's $\varphi(\mu):=f(\mu)-f_0(\mu)$, $f_0(\mu)$ and $f(\mu)$ admit
meromorphic continuations  $\varphi(z)$, $f_0(z)$ and $f(z)$ on $\Omega_+$, such that
\begin{equation}\label{eq:3.16}
M_\pi(f-f_0, \Omega_+)=\kappa.
\end{equation}
Let $s=b_1fb_2$ and $K=b_1f_0b_2$.  Then the equality~\eqref{eq:3.16} yields
\[
M_\pi(s-K, \Omega_+)\leq\kappa.
\]
Since $K\in H_\infty^{\ptq}$, then
\[
\kappa':=M_\pi(s, \Omega_+)=M_\pi(s-K, \Omega_+)\leq\kappa.
\]
Taking into account that $\|s\|_\infty=\|f\|_\infty\leq 1$, one obtains $s\in\mathcal{S}_{\kappa'}$.
Moreover,  the condition~\eqref{eq:3.16} is equivalent to the condition~\eqref{eq:3.15},
i.e. $s\in\mathcal{TS}_\kappa(b_1, b_2, K)$.

Conversely, if  $s\in\mathcal{S}_{\kappa'}^{p\times q}$ with $\kappa'\leq \kappa$ and the
condition \eqref{eq:3.15} is in force, then for $f=b^{-1}_1sb_2^{-1}$, $f_0=b_1^{-1}Kb_2^{-1}$ one obtains
that \eqref{eq:3.16} holds
and $\|f\|_\infty\leq 1$. Therefore, $f\in\mathcal{N}_\kappa(f_0)$.
\end{proof}

\begin{lemma}
\label{lem:oct28a11}
Let $f_0\in L^{p\times q}_\infty$, $\Gamma=\Gamma(f_0)$, $\{b_1, b_2\}\in den(f_0)$, $K=b_1f_0b_2$ and
let ${\bf P}$ be the operator in $\cH(b_1)\oplus\cH_*(b_2)$, defined by formulas \eqref{eq:I26} and
\eqref{eq:I30}. Then
\[
\label{eq:nu_P} \nu_-({\bf P})=\nu_-(I-\Gamma^*\Gamma).
\]
Moreover, if $\nu_-(I-\Gamma^*\Gamma)<\infty$, then
\[
\label{eq:rho_P}
0\in\rho({\bf P})\Longleftrightarrow 0\in\rho(I-\Gamma^*\Gamma).
\]
\end{lemma}
\begin{proof}
Let us decompose the spaces $H_2^q$ and $(H_2^p)^\perp$:
\[
H_2^q=b_2(H_2^q)\oplus\cH(b_2),\quad
(H_2^p)^\perp=\cH_*(b_1)\oplus b_1(H_2^p)^\perp
\]
and let us decompose the operator $\Gamma:H_2^q\to(H_2^p)^\perp$, accordingly
\begin{equation}
\label{eq:oct28a11}
\Gamma\stackrel{def}{=}\left(
                                      \begin{array}{cc}
                                        \Gamma_{11} & \Gamma_{12} \\
                                        0            & \Gamma_{22} \\
                                      \end{array}
                                    \right):\,
\begin{array}{c}
                                     b_2(H_2^q) \\
                                     \oplus \\
                                     \cH(b_2)
                                   \end{array}\to\begin{array}{c}
                                                   \cH_*(b_1) \\
                                                   \oplus \\
                                                  b_1^*(H_2^p)^\perp
                                                 \end{array},
\end{equation}
where the operators
\[
\Gamma_{11}: b_2(H_2^q)\to\cH_*(b_1),\quad  \Gamma_{12}:
\cH(b_2)\to\cH_*(b_1),\quad \Gamma_{22}: \,\cH(b_2)\to b_1^*(H_2^p)^\perp
\]
are defined by the formulas
\begin{equation}\label{eq:I26G}
\begin{split}
\Gamma_{11}h_+&=\Pi_{\cH_*(b_1)}Kh_+,\qquad h_+\in b_2(H_2^q),\\
\Gamma_{12}h_2&=\Pi_{\cH_*(b_1)}Kh_2,\,\qquad h_2\in\cH(b_2),\\
\Gamma_{22}h_2&=(b_1^*\Pi_-b_1)Kh_2,\,\,\,\quad h_2\in\cH(b_2).
\end{split}
\end{equation}
It follows from~\eqref{eq:oct28a11}, \eqref{eq:I26G} and~\eqref{eq:I26} that the operator
$\Gamma:H_2^q\to(H_2^p)^\perp$ and the operator
\[
{\bf K}\stackrel{def}{=}\left(
                                      \begin{array}{cc}
                                        K_{11} & K_{12} \\
                                        0            & K_{22} \\
                                      \end{array}
                                    \right):\,
\begin{array}{c}
                                     H_2^q \\
                                     \oplus \\
                                     \cH_*(b_2)
                                   \end{array}\to\begin{array}{c}
                                                   \cH(b_1) \\
                                                   \oplus \\
                                                  (H_2^p)^\perp
                                                 \end{array}
\]
are connected by
\[
\Gamma=\cM_{b_1^*}|_{b_1(H_2^p)^\perp}{\bf K}\cM_{b_2}|_{H_2^q}.
\]
and, hence, the operators $\Gamma$ and ${\bf K}$ are unitary equivalent. Now the statements are implied
by~\cite[Lemma~5.10]{DD14}.
\end{proof}

\begin{theorem}\label{thm:4.3}
Let $f_0\in L^{p\times q}_\infty$, $\Gamma=\Gamma(f_0)$, $\kappa\in\mathbb{Z}_+$, $\{b_1, b_2\}\in den(f_0)$,
$K=b_1f_0b_2$, let ${\bf P}$ be defined by formulas \eqref{eq:I30},  let (H1)--(H4) be in force, let  the mvf
$W(z)$ be defined by~\eqref{eq:2.11}
 and let
\begin{equation}\label{eq:4.39G}
\mathfrak{A}(\mu)=\begin{bmatrix} b_1(\mu)^{-1}&0\\0&b_2(\mu)   \end{bmatrix}
W(\mu),
\quad\mu\in\gh_{b_1}\cap\gh_{b_2^\#}\cap\Omega_0.
 \end{equation}
Then:
\begin{enumerate}
  \item [(1)] $\mathfrak{A}\in\mathfrak{M}^r_\kappa(j_{pq})$;
  \item [(2)] $\cN_\kappa(f_0)\ne\emptyset$ if and only if $\kappa\ge\kappa_1:=\nu_-(I-\Gamma^*\Gamma)$;
  \item [(3)] $\cN_\kappa(f_0)=T_{\mathfrak{A}}[\cS_{\kappa-\kappa_1}^{\ptq}]$,
    \item [(4)] $\cN\cT_\kappa(f_0)=\cup_{k=\kappa_1}^\kappa T_{\mathfrak{A}}[\cS_{k-\kappa_1}^{\ptq}]$.
\end{enumerate}
\end{theorem}
\begin{proof}
(1)
By~\cite[Theorem 4.2]{DD14} the rows of $W(z)$ admit  factorizations
\[
\begin{bmatrix} w_{11}& w_{12}\end{bmatrix}=b_1\begin{bmatrix} a_{11}& a_{12}\end{bmatrix},
\]
\[
\begin{bmatrix} w_{21}& w_{22}\end{bmatrix}=b_2^{-1}\begin{bmatrix} a_{21}& a_{22}\end{bmatrix},
\]
where $a_{11}\in(H_2^{p\times p})^\perp, a_{12}\in(H_2^{p\times q})^\perp$, $a_{21}\in H_2^{q\times p}$, $a_{22}\in H_2^{q\times q}$ and
\[
s_{21}=-w_{22}^{-1}w_{21}=-a_{22}^{-1}a_{21}\in\mathcal{S}_\kappa^{p\times q}.
\]
If the mvf's $b_\ell^{-1}$, $s_\ell$, $b_r$, $s_r$ are determined by Krein-Langer factorizations of $s_{21}$
\[
s_{21}=b_\ell^{-1}s_\ell=s_rb_r^{-1},
\]
then in accordance with~\cite[Theorem 4.3]{DD14} (see (4.26), (4.27))
\[
a_2:=b_\ell a_{22}^{-1}\in \mathcal{S}_{out}^{q\times q}, \quad
a_1:=(a_{11}^\#)^{-1}b_r\in \mathcal{S}_{out}^{p\times p}.
\]
Thus
\[
\mathfrak{A}(z)=\begin{bmatrix}b_1^{-1}&0\\0&b_2\end{bmatrix}
\begin{bmatrix}w_{11}&w_{12}\\w_{21}&w_{22}\end{bmatrix}
\]
belongs to the class $\mathfrak{M}^r_\kappa(j_{pq})$.

(2)
 By Theorem~\ref{thm:3.1}  $\mathcal{N}_\kappa(f_0)$ is nonempty if and only if  $\mathcal{TS}_\kappa(b_1,b_2,K)$
 is nonempty. Therefore (2) is implied by Theorem~\ref{thm:nov11a13} and Lemma~\ref{lem:oct28a11}.

(3) The statement (3) follows from the formula~\eqref{eq:ST_descr} proved in Theorem~\ref{thm:nov11a13} and
from the equivalence
 \[
     f\in \mathcal{N}_\kappa(f_0)\Longleftrightarrow
     b_1 fb_2\in \mathcal{TS}_\kappa(b_1,b_2,K)
     T_{W}[\cS_{\kappa-\kappa_1}]
 \]
(Theorem~\ref{thm:3.1}). This means that for every $f\in \mathcal{N}_\kappa(f_0)$ the mvf $s= b_1 fb_2$  belongs
to $\mathcal{TS}_\kappa(b_1,b_2,K)$ and hence it admits the representation
\[
    s=(w_{11}\varepsilon+w_{12})(w_{21}\varepsilon+w_{22})^{-1}=T_{W}[\varepsilon]
\]
for some $\varepsilon\in\cS_{\kappa-\kappa_1}$. Therefore, the mvf $f=b_1^{-1}sb_2^{-1}$ can be represented as
\[
    f=b_1^{-1}(w_{11}\varepsilon+w_{12})(b_2w_{21}\varepsilon+b_2w_{22})^{-1}=T_{\mathfrak A}[\varepsilon].
\]

(4) As follows from (2) $\cN_{\kappa'}(f_0)=\emptyset$ for $\kappa'<\kappa_1$.
Therefore, (4) is implied by~\eqref{eq:3.13} and by the statement~(3).
\end{proof}

\section{Resolvent matrix in the case of a rational mvf $f_0$}
Assume now that $\Omega_+=\dD$ and $f_0$ is a rational mvf with a minimal
realization
\begin{equation}\label{eq:4.0}
f_0(z)=C(zI_n-A)^{-1}B,
\end{equation}
where $n\in\mathbb{N}$, $A\in\mathbb{C}^{n\times n}$, $B\in\mathbb{C}^{n\times q}$, $C\in\mathbb{C}^{p\times n},$ �
\begin{equation}\label{eq:4.1}
 \quad \sigma(A)\subset\mathbb{D}.
\end{equation}
Then the corresponding Hankel operator $\Gamma=\Gamma(f_0):H_2^q\to (H_2^p)^\perp$ in~\eqref{eq:3.1}
has the following matrix representation $(\gamma_{j+k-1})_{j,k=1}^\infty$ in the standard
basises $\{e^{ijt}\}_{j=0}^\infty$ and $\{e^{-ikt}\}_{k=1}^\infty$:
\[
    (\gamma_{j+k-1})_{j,k=1}^\infty=(CA^{j+k-2}B)_{j,k=1}^\infty=\Omega\Xi,
\]
where $\gamma_j$ are given by~\eqref{eq:gf_1}.
\begin{equation}\label{eq:4.2C}
\Xi=\begin{bmatrix}B& AB& \ldots&A^{n-1}B\end{bmatrix} \quad
\mbox{and}\quad
\Omega=\begin{bmatrix}CA^0\\\vdots\\CA^{n-1}\end{bmatrix}.
\end{equation}

Representation \eqref{eq:4.0} is called minimal, if the dimension of
the matrix $A$ in~\eqref{eq:4.0} is minimal. As is known
see~\cite[Thm 4.14]{BGR90} the representation \eqref{eq:4.0} is
minimal if and only if the pair $(A,B)$ is controllable and the pair
$(C,A)$ is observable, i.e.
\begin{equation}\label{eq:4.2}
\ran\Xi=\mathbb{C}^n \quad \mbox{and}\quad \ker\Omega=\{0\},
\end{equation}
The  controllability gramian  $P$  and the observability gramian
$Q$, defined by
\[
P=\sum_{k=0}^\infty A^kBB^*(A^*)^k=\Xi\Xi^*, \quad
Q=\sum_{k=0}^\infty (A^*)^kCC^*(A)^k=\Omega^*\Omega,
\]
are solutions to the following Lyapunov-Stein equations
\begin{equation}\label{eq:Gram}
P-APA^*=BB^*, \quad Q-A^*QA=C^*C.
\end{equation}

As was shown in~\cite[Remark 4.2]{DD10} a denominator of the mvf $f_0(z)$ may be selected as $(I_p, b_2)$, where
\begin{equation}\label{eq:4.3}
b_2(z)=I_q-(1-z)B^*(I_n-zA^*)^{-1}P^{-1}(I_n-A)^{-1}B
\end{equation}
Straightforward calculations show that
\begin{equation}\label{eq:4.5}
(zI_n-A)^{-1}Bb_2(z)=P(I_n-A^*)(I_n-zA^*)^{-1}P^{-1}(I_n-A)^{-1}B.
\end{equation}
Since the mvf $b_2(z)$ is inner, then $b_2(z)^{-1}=b_2(\frac{1}{\bar{z}})^*$, and hence
\begin{equation}\label{eq:4.6}
b_2(z)^{-1}=I_q+(1-z)B^*(I_n-A^*)^{-1}P^{-1}(zI_n-A)^{-1}B.
\end{equation}
\begin{proposition}\label{thm:4.1}
Let $f_0(z)$ be a mvf of the form \eqref{eq:4.0}, where
$A\in\mathbb{C}^{n\times n}$, $B\in\mathbb{C}^{n\times q}$,
$C\in\mathbb{C}^{p\times n}$ satisfy~\eqref{eq:4.1}
and~\eqref{eq:4.2} and let
 \begin{equation}\label{eq:4.8A}
 M=\begin{bmatrix}-A & 0\\ 0&I_n\end{bmatrix},\quad
  N=\begin{bmatrix}-I_n & 0\\ 0&A^*\end{bmatrix},\quad
 \Lambda=\begin{bmatrix}-Q & I_n\\ I_n&-P\end{bmatrix},
 \end{equation}
  \begin{equation}\label{eq:4.8B}
  G(z)=\begin{bmatrix}C & 0\\ 0& B^*\end{bmatrix}(M-zN)^{-1}.
 \end{equation}
Assume that $1\not\in\sigma(PQ)$. Then:
\begin{enumerate}
  \item [(1)] $\cN_\kappa(f_0)\ne\emptyset$ if and only if $\kappa_1:=\nu_-(I-PQ)\le\kappa$;
  \item [(2)] If (1) holds then the matrix $\Lambda$ is invertible and
 $\mathcal{N}_\kappa(f_0)=T_{\mathfrak{A}}[\cS_{\kappa-\kappa_1}]$,
 where
 \begin{equation}\label{eq:4.8}
 \mathfrak{A}(\mu)=I_m-(1-\mu)G(\mu)
\Lambda^{-1}G(1)^*j_{pq};
 \end{equation}
  \item [(3)]
The mvf $\mathfrak{A}(\mu)$ is a generalized right
$\gamma$-generating mvf of the class
$\mathfrak{M}_{\kappa_1}^r(j_{pq})$.
\end{enumerate}
\end{proposition}
The statements (1), (2) of Proposition~\ref{thm:4.1} and the formula
\eqref{eq:4.8} for the resolvent matrix $\mathfrak{A}(\mu)$ are well
known from~\cite[Theorem~20.5.1]{BGR90}. We will deduce now the
formula~\eqref{eq:4.8} from the general formula~\eqref{eq:2.11} for
the resolvent matrix of the problem $\mathbf{TSP}_\kappa(I_p,b_2,K)$
with
\begin{equation}\label{eq:4.9}
 K(z)=f_0(z)b_2(z)=C(zI_n-A)^{-1}Bb_2(z).
 \end{equation}
\begin{proof}
(1) By Theorem~\ref{thm:3.1} $f\in \mathcal{N}_\kappa(f_0)$ if and
only if $s=fb_2\in \mathcal{TS}_\kappa(I_p,b_2,K)$. Alongside with
the problem $\mathbf{TSP}_\kappa(I_p,b_2,K)$ consider the problem
$\mathbf{GSTP}_\kappa(I_p,b_2,K)$ (see Remark~\ref{rem:GSTP}). As is
known~\cite[Theorem 5.17]{DD10} these problems have the same
resolvent matrix. Assume that $s\in\cS_\kappa(I_p,b_2,K)$,
see~\eqref{eq:Sdef}. The conditions $s\in\cS_\kappa^{\ptq}$ and
$(s-K)b_2^{-1}\in\cN_{+,\kappa}$ are equivalent to the equalities
\[
M_\pi(s,\Omega_+)=M_\pi((s-K)b_2^{-1},\Omega_+)=\kappa.
\]
By the cancellation lemma from~\cite[Lemma 5.5]{DD10}
 \begin{equation}\label{eq:4.10}
M_\pi((s_\ell-b_\ell K)b_2^{-1},\Omega_+)=M_\pi(s_\ell,\Omega_+)=0.
 \end{equation}
By~\eqref{eq:4.6} and~\eqref{eq:4.9} the expression $(s_\ell-b_\ell
K)b_2^{-1}$ takes the form
 \[
     s_\ell(I_q+(1-z)B^*(I_n-A^*)^{-1}P^{-1}(zI_n-A)^{-1}B)-b_\ell C(zI_n-A)^{-1}B.
 \]
Since $ (1-z)(zI_n-A)^{-1}=-I_n+(I_n-A)(zI_n-A)^{-1} $, the
condition \eqref{eq:4.10} can be rewritten as
 \begin{equation}\label{eq:4.12}
\{s_\ell B^*(I_n-A^*)^{-1}P^{-1}(I_n-A)-b_\ell C\}(zI_n-A)^{-1}B\in\mathcal{N}_+.
 \end{equation}
 Since the pair $(A, B)$ is controllable, then \eqref{eq:4.12} is equivalent to
  \begin{equation}\label{eq:4.13}
\begin{bmatrix} b_\ell & -s_\ell\end{bmatrix}
\begin{bmatrix} C \\ B^*(I_n-A^*)^{-1} P^{-1}(I_n-A)\end{bmatrix}(zI_n-A)^{-1}\in\mathcal{N}_+.
 \end{equation}
 Thus, the condition \eqref{eq:4.13} can be rewritten as
  \begin{equation}\label{eq:4.14}
\begin{bmatrix} b_\ell & -s_\ell\end{bmatrix} F\in\mathcal{N}_+,
 \end{equation}
where
  \begin{equation}\label{eq:4.15}
F(z)=\widetilde{C}(A-zI_n)^{-1}, \quad \widetilde{C}=\begin{bmatrix} C \\ B^*(I_n-A^*)^{-1}P^{-1}(I_n-A)\end{bmatrix}.
 \end{equation}
Thus, the problem $\mathbf{GSTP}_\kappa(I_p,b_2,K)$ is equivalent to the one-sided
interpolation problem~\eqref{eq:4.14} considered in~\cite{DD10}.
As was shown in~\cite[(1.14)]{DD10} the Pick matrix $\wt P$, corresponding to the problem~\eqref{eq:4.14}
is the unique solution of the Lyapunov-Stein equation
 \begin{equation}\label{eq:4.19}
A^*\wt PA-\wt P=\widetilde{C}^*j_{pq}\widetilde{C}
 \end{equation}
and the problem~\eqref{eq:4.14} is solvable if and only if
\[
    \kappa_1:=\nu_-(\wt P)\le\kappa.
\]
Since by~\eqref{eq:Gram}
\begin{equation}\label{eq:4.22}
\begin{split}
\widetilde{C}^*j_{pq}\widetilde{C}&=
C^*C-(I_n-A^*)P^{-1}(I_n-A)^{-1}BB^*(I_n-A^*)^{-1}P^{-1}(I_n-A)
\\&=C^*C-(I_n-A^*)P^{-1}-A^*P^{-1}(I_n-A)
\\&=(Q-P^{-1})-A^*(Q-P^{-1})A,
\end{split}
 \end{equation}
then
 \begin{equation}\label{eq:4.21}
\wt P=P^{-1}-Q=P^{-1/2}(I-P^{1/2}QP^{1/2})P^{-1/2}.
 \end{equation}
Notice, that in \eqref{eq:4.22} we use the equality
\[
-(I_n-A)^{-1}BB^*(I_n-A^*)^{-1}=-(I_n-A)^{-1}P-PA^*(I_n-A^*)^{-1},
 \] 
It follows from~\eqref{eq:4.21} and Theorem~\ref{thm:nov11a13} that
$\mathcal{TS}_\kappa(I_p,b_2,K)\ne\emptyset$ if and only if
\[
    \kappa_1:=\nu_-(I-P^{1/2}QP^{1/2})\le\kappa.
\]
Now it remains to note that
$\sigma(I-P^{1/2}QP^{1/2})=\sigma(I-PQ)$. In view of
Theorem~\ref{thm:3.1} this proves (1).

(2)
By \cite[Theorem~3.1 and Theorem~5.17]{DD10} the resolvent matrix
$\wt W(z)$, which describes the set
$\mathcal{TS}_\kappa(I_p,b_2,K_0)$ via the
formula~\eqref{eq:STdescr} takes the form
  \[
    \wt W(z)=I_m-(1-z)F(z)\wt P^{-1}F(1)^*j_{pq},
 \]
where $\wt P$ is given by \eqref{eq:4.19}. As was shown
in~\cite[Lemma~4.8]{DD10} $\wt W\in\cU_\kappa^r(j_{pq})$. Let us set
\begin{equation}\label{eq:6.1}
    \wt {\mathfrak A}(\mu)
    :=  \begin{bmatrix}
            I_p&0\\0&b_2(\mu)
        \end{bmatrix}
        W(\mu)
\end{equation}
and show that the mvf
\begin{equation}\label{eq:6.2}
    \wt {\mathfrak A}(\mu)
    :=  \begin{bmatrix}
            I_p&0\\0&b_2(\mu)
        \end{bmatrix}
               +(\mu-1)\begin{bmatrix}
            I_p&0\\0&b_2(\mu)
        \end{bmatrix}
        F(\mu)\wt P^{-1}F(1)^*j_{pq}
\end{equation}
coincides with the mvf  ${\mathfrak A}$ from~\eqref{eq:4.8}.
It follows from~\eqref{eq:4.5} that
 \[
 b_2(\mu)^{-1}B^*(I_n-\mu A^*)^{-1}=B^*(I_n-A^*)^{-1}P^{-1}(\mu I_n-A)^{-1}(I_n-A)P.
 \]
and hence
\[ 
 b_2(\mu)B^*(I_n-A^*)^{-1}P^{-1}(\mu I_n-A)^{-1}(I_n-A)=B^*(I_n-\mu A^*)^{-1}P^{-1}.
\]
In view of 
\eqref{eq:4.3},  \eqref{eq:4.15}, \eqref{eq:4.8A} and \eqref{eq:4.8B} this implies
\begin{equation}\label{eq:6.3}
\begin{bmatrix}
            I_p&0\\0&b_2(\mu)
        \end{bmatrix}
        F(\mu)=
    \begin{bmatrix}
        C(\mu I_n-A)^{-1}\\ B^*(I_n-\mu A^*)^{-1}P^{-1}
    \end{bmatrix}
    =G(\mu)
    \begin{bmatrix}
        I_p\\ P^{-1}
    \end{bmatrix}.
   \end{equation}
Next, in view of 
\eqref{eq:4.15}  and \eqref{eq:4.3}
\begin{equation}\label{eq:F1*}
    F(1)^*=
    \begin{bmatrix}
        (I_n-A^*)^{-1}C^*& P^{-1}(I_n-A)^{-1}B
    \end{bmatrix}
    =\begin{bmatrix}
        I_n& P^{-1}
    \end{bmatrix}
    G(1)^*,
\end{equation}
\begin{equation}\label{eq:Ib_2}
\begin{split}
    \begin{bmatrix}
            I_p&0\\0&b_2(\mu)
    \end{bmatrix}&=I_m-(1-\mu)
    \begin{bmatrix} 0& 0
        \\ 0& B^*(I_n-zA^*)^{-1}P^{-1}(I_n-A)^{-1}B
    \end{bmatrix}\\
        &=I_m-(1-\mu)G(\mu)
    \begin{bmatrix} 0& 0
                \\ 0 & -P^{-1}
    \end{bmatrix}
    G(1)^*j_{pq}.
    \end{split}
\end{equation}
Substituting~\eqref{eq:6.3},~\eqref{eq:F1*} and~\eqref{eq:Ib_2} into~\eqref{eq:6.2} one obtains
\[
    \wt {\mathfrak A}(\mu)=I_m-(1-\mu)G(\mu)
    \begin{bmatrix} \wt P^{-1}& \wt P^{-1}P^{-1}
                \\ P^{-1}\wt P^{-1}  & -P^{-1}+P^{-1}\wt P^{-1}P^{-1}
    \end{bmatrix}
    G(1)^*j_{pq}.
\]
In view of the equality
\[
 \begin{bmatrix} \wt P^{-1}& \wt P^{-1}P^{-1}
                \\ P^{-1}\wt P^{-1}  & -P^{-1}+P^{-1}\wt P^{-1}P^{-1}
    \end{bmatrix}=\Lambda^{-1}
 \]
 this proves the formula~\eqref{eq:4.8}.

By \cite[Theorem~3.1 and Theorem~5.17]{DD10} and
Theorem~\ref{thm:nov11a13} the set
$\mathcal{TS}_\kappa(I_p,b_2,K_0)$ is described by the formula
\[
\cT\cS_\kappa(b_1,b_2;K)
=T_W[{\cS}_{\kappa-\kappa_1}^{\ptq}]
=\{T_W[\varepsilon]:\,\varepsilon\in\cS_{\kappa-\kappa_1}^{\ptq}\}.
\]
Therefore, the statement (3) is implied by Theorem~\ref{thm:4.3} (3).

(3) By Theorem~\ref{thm:4.3} the mvf $\wt {\mathfrak A}(\mu)$
defined by \eqref{eq:6.1} belongs to the class
$\mathfrak{M}_{\kappa_1}^r(j_{pq})$.
\end{proof}

\end{document}